\documentstyle{amsppt}

\refstyle{A}

\nologo

\hoffset .25 true in
\voffset .2 true in

\hsize=6.1 true in
\vsize=8.5 true in

\define\Adot{\bold A^\bullet}
\define\Bdot{\bold B^\bullet}
\define\Fdot{\bold F^\bullet}
\define\map{{f\circ\pi_1+g\circ\pi_2}}
\define\lotimes{\ {\overset L\to \otimes}\ }
\define\lboxtimes{\ {\overset L\to \boxtimes}\ }
\define\piten{{\pi_1^*\Adot\lotimes\pi_2^*\Bdot}}
\define\rhom{R\bold{Hom}^\bullet}

\topmatter

\title The Sebastiani-Thom Isomorphism in the Derived category \endtitle

\author David B. Massey \endauthor

\address{David B. Massey, Dept. of Mathematics, Northeastern University, Boston, MA, 02115, USA} \endaddress

\email{DMASSEY\@NEU.edu}\endemail



\endtopmatter

\document

\vskip .2in

\noindent\S0. {\bf Introduction}  

\vskip .2in

Let $f:X\rightarrow\Bbb C$ and $g:Y\rightarrow\Bbb C$ be analytic functions. Let $\pi_1$ and $\pi_2$ denote the projections
of $X\times Y$ onto $X$ and $Y$, respectively.  In [{\bf S-T}], Sebastiani and Thom  prove that the cohomology of the Milnor
fibre of $\map$ is isomorphic to the tensor product of the cohomologies of the Milnor fibres of $f$ and $g$ (with a shift in
degrees); they prove this in the case where
$X$ and $Y$ are smooth and $f$ and $g$ have isolated critical points. In addition, they prove that the monodromy isomorphism
induced by $\map$ is the tensor product of those induced by $f$ and $g$. The point, of course, is to break up the 
complicated critical activity of $\map$ into more manageable pieces. Sebastiani-Thom-type results have been
proved by N\'emethi [{\bf N1}], [{\bf N2}], Oka [{\bf O}], and Sakamoto [{\bf S}].

\vskip .1in

In this paper, we prove that this Sebastiani-Thom Isomorphism exists regardless of how singular the spaces $X$ and $Y$ may
be, regardless of how bad the critical loci of $f$ and $g$ are, and regardless of what coefficients one uses. Moreover, we
prove that the Sebastiani-Thom Isomorphism is actually a consequence of a natural isomorphism in the derived category of
bounded, constructible complexes of sheaves on $X\times Y$. 

\vskip .1in

To state our result precisely, we must introduce more notation -- all of which can be found in [{\bf K-S}].

\vskip .1in

 Let $R$ be a regular Noetherian ring with finite Krull dimension (e.g., $\Bbb Z,
\Bbb Q, \text{or}\ \Bbb C$). Let $\Adot$ and $\Bdot$ be bounded, constructible complexes of sheaves of $R$-modules on
$X$ and $Y$, respectively. Recall that, in this situation, $\Adot\lboxtimes\Bdot := \piten$.

\vskip .1in

Let  $p_1$ and $p_2$ denote the projections of $V(f)\times V(g)$ onto $V(f)$ and $V(g)$, respectively, and let
$k$ denote the inclusion of $V(f)\times V(g)$ into $V(f\circ\pi_1+g\circ\pi_2)$.

\vskip .1in

If $h:Z\rightarrow\Bbb C$ is an analytic function, and $\Fdot$ is a complex on $Z$, then $\phi_{{}_h}\Fdot$ denotes the sheaf
of vanishing cycles of $\Fdot$ along $h$. Here, $\phi_{{}_h}\Fdot$ is defined as in 8.6.2 of [{\bf K-S}] and, hence, is shifted
by
$1$ from the more traditional definitions, i.e., in this paper, the stalk cohomology of
$\phi_{{}_h}\Fdot$ in degree $i$ is the degree $i$ relative hypercohomology of a small ball modulo the Milnor fibre with
coefficients in $\Fdot$. Thus, in the constant $\Bbb Z$-coefficient case, $H^i(\phi_{{}_h}{\Bbb Z}^\bullet_X)_{\bold x}\cong
\widetilde H^{i-1}(F_{{}_{h, \bold x}})$, where $\widetilde H$ denotes reduced cohomology and
$F_{{}_{h, \bold x}}$ denotes the Milnor fibre of $h$ at $\bold x$.

\vskip .2in

We prove:

\vskip .1in

\noindent{\bf Theorem (Sebastiani-Thom Isomorphism)}. {\it There is a natural
isomorphism
$$
k^*\phi_{{}_\map}\big(\Adot\lboxtimes\Bdot\big)\ \cong\ \phi_{{}_f}\Adot\lboxtimes \phi_{{}_g}\Bdot,
$$
and this isomorphism commutes with the corresponding monodromies. 

Moreover, if we let $\bold p:=(\bold x, \bold y)\in
X\times Y$ be such that
$f(\bold x) = 0$ and $g(\bold y)=0$, then, in an open neighborhood of $\bold p$, the
complex 
$\phi_{{}_\map}\big(\piten\big)$ has support contained in $V(f)\times V(g)$, and, in any open set in which we have this
containment, there are natural isomorphisms
$$
\phi_{{}_\map}\big(\Adot\lboxtimes\Bdot\big)\ \cong\ k_!(\phi_{{}_f}\Adot\lboxtimes \phi_{{}_g}\Bdot)\ \cong\
k_*(\phi_{{}_f}\Adot\lboxtimes \phi_{{}_g}\Bdot).
$$}

\vskip .2in

\noindent \S1. {\bf Proof of the General Sebastiani-Thom Isomorphism}

\vskip .2in

The proof of the Sebastiani-Thom Isomorphism uses two Morse-theoretic lemmas and two formal, derived category propositions.
First, however, we need to discuss the definition of the vanishing cycles that we shall use.

Kashiwara and Schapira define the vanishing cycles in 8.6.2 of [{\bf K-S}]. However, we shall use the more down-to-Earth
characterization (via natural equivalence) given in Exercise VIII.13 of [{\bf K-S}]. Hence, we use as our definition:
$\phi_{{}_h}\Fdot = \big(R\Gamma_{{}_{\{\operatorname{Re} h\leqslant 0\}}}\Fdot\big)_{|_{V(h)}}$. (We have reversed the
inequality used in [{\bf K-S}]. We do this for aesthetics only -- we prefer to think of the vanishing cycles in terms of a
ball modulo the Milnor fibre over a small {\bf positive} value of the function.) The monodromy isomorphism is easy to
describe: for all $\theta$, there is an isomorphism $\big(R\Gamma_{{}_{\{\operatorname{Re} h\leqslant
0\}}}\Fdot\big)_{|_{V(h)}}\cong \big(R\Gamma_{{}_{\{\operatorname{Re} e^{-i\theta}h\leqslant 0\}}}\Fdot\big)_{|_{V(h)}}$,
and the monodromy isomorphism results when $\theta = 2\pi$.

\vskip .3in

We continue with our notation from the introduction.

\vskip .5in

\noindent{\bf Lemma 1.1}. {\it Let $\Cal S$ (resp. $\Cal S^\prime$) denote a complex Whitney stratification of $X$ (resp.
$Y$) with respect to which $\Adot$ (resp. $\Bdot$) is constructible. Let $\Sigma_{{}_{\Cal S}}f$ and $\Sigma_{{}_{\Cal
S^\prime}}g$ denote the stratified critical loci. Then, 

$$
\operatorname{supp} \phi_{{}_\map}\big(\Adot\lboxtimes\Bdot\big)\ \subseteq\ \big(\Sigma_{{}_{\Cal S}}f\times \Sigma_{{}_{\Cal
S^\prime}}g\big)\cap V(\map);
$$
in particular, if $\Sigma_{{}_{\Cal S^\prime}}g\subseteq V(g)$, then 
$$
\operatorname{supp} \phi_{{}_\map}\big(\Adot\lboxtimes\Bdot\big)\ \subseteq\ V(f)\times V(g).
$$

\vskip .1in

Moreover, if $\bold p:=(\bold x, \bold y)\in
X\times Y$ is such that
$f(\bold x) = 0$ and $g(\bold y)=0$, then, near
$\bold p$,
$$
\operatorname{supp} \phi_{{}_\map}\big(\Adot\lboxtimes\Bdot\big)\ \subseteq\ V(f)\times V(g).
$$ }

\vskip .2in

\noindent{\it Proof}. As the complexes $\Adot$ and $\Bdot$ are constructible with respect to $\Cal S$ and $\Cal S^\prime$,
$\Adot\lboxtimes\Bdot$ is constructible with respect to the product stratification, $\Cal S\times\Cal S^\prime$. The support of
$\phi_{{}_\map}\big(\Adot\lboxtimes\Bdot\big)$ is contained in the stratified critical locus of $\map$, which is trivially
seen to be equal to the product of the stratified critical loci of
$f$ and $g$. Finally, near $\bold x$ and $\bold y$, these stratified critical loci are contained in $V(f)$ and
$V(g)$, respectively.\qed

\vskip .3in

Recall that $k:V(f)\times V(g)\hookrightarrow V(f\circ\pi_1+g\circ\pi_2)$ denotes the inclusion. Let
$q:V(f\circ\pi_1+g\circ\pi_2)\hookrightarrow X\times Y$ and $m: V(f)\times V(g)\hookrightarrow X\times Y$ also denote the
inclusions, so that $m= q\circ k$.

\vskip .5in

\noindent{\bf Lemma 1.2}. {\it Let $P:=\{\bold x\in X\ |\ \operatorname{Re} f(\bold x)\leqslant 0\}$, let
$Q:=\{\bold y\in Y\ |\ \operatorname{Re} g(\bold y)\leqslant 0\}$, and let $Z:=\{(\bold x, \bold y)\in X\times Y\ |\
 \operatorname{Re}(f(\bold x)+g(\bold y))\leqslant 0\}$. Note that $P\times Q\subseteq Z$.

The natural map
$$ R\Gamma_{{}_{P\times Q}}(\Adot\lboxtimes\Bdot) \ \rightarrow\ R\Gamma_{{}_{Z}}(\Adot\lboxtimes\Bdot)
$$ induces a natural isomorphism
$$ m^*R\Gamma_{{}_{P\times Q}}(\Adot\lboxtimes\Bdot) \ \cong\ m^*R\Gamma_{{}_{Z}}(\Adot\lboxtimes\Bdot)\ \cong\
k^*\phi_{{}_\map}\big(\Adot\lboxtimes\Bdot\big).
$$

}

\vskip .2in

\noindent{\it Proof}. As in the first lemma, we use nothing about $\Adot\lboxtimes\Bdot$ other than the fact that it is
constructible with respect to the product stratification; let us use $\Fdot$ to denote $\Adot\lboxtimes\Bdot$. 

From the definition of the vanishing cycles, we have 
$$m^*R\Gamma_{{}_{Z}}(\Fdot)\ =\ (q\circ k)^*R\Gamma_{{}_{Z}}(\Fdot)\ \cong\ k^*q^*R\Gamma_{{}_{Z}}(\Fdot)\ 
\cong\ k^*\phi_{{}_\map}\big(\Fdot\big).$$

Let $\bold p\in V(f)\times V(g)$. We need to show that we have the
isomorphism 
$$H^*(R\Gamma_{{}_{P\times Q}}(\Fdot))_{\bold p} @>\cong>>H^*(R\Gamma_{{}_{Z}}(\Fdot))_{\bold p}.$$
 
\vskip .1in

Let $\Theta:X\times Y\rightarrow \Bbb R^2$ be given by $\Theta(\bold x, \bold y) := (\operatorname{Re} f(\bold x),
\operatorname{Re} g(\bold y))$. Use $u$ and $v$ for the coordinates in $\Bbb R^2$. Let $C:=\{(u,v)\ |\ u\leqslant 0,
v\leqslant 0\}$, and
$D:=\{(u,v)\ |\ u+v\leqslant 0\}$. What we need to show is that we have an isomorphism 
$$H^*(R\Gamma_{{}_{\Theta^{-1}(C)}}(\Fdot))_{\bold p}@>\cong>>H^*(R\Gamma_{{}_{\Theta^{-1}(D)}}(\Fdot))_{\bold
p}. \tag{$\dagger$}$$

In a small neighborhood of $\bold p$, the map
$\Theta$ is a stratified submersion over the complement of $\{(u,v)\ |\ uv=0\}$ (the coordinate-axes); for a critical
point of $\operatorname{Re} f$ (resp. $\operatorname{Re} g$) on a stratum occurs at a critical point of $f$ (resp. $g$) on
the stratum. The desired result will follow by {\it moving the wall} (see [{\bf G-M}]); essentially, one deforms the region
$D$ to the region
$C$ by the obvious homotopy, and lifts this deformation up to $X\times Y$ via $\Theta$. 

\vskip .1in

To avoid the critical values along the $u$ and $v$ axes, it is slightly easier to work with the complements of $C$ and $D$.
Note that, since we have the morphism of distinguished triangles
$$
\rightarrow R\Gamma_{{}_{\Theta^{-1}(C)}}(\Fdot)\rightarrow\Fdot\rightarrow R\Gamma_{{}_{\Theta^{-1}(\Bbb
R^2-C)}}(\Fdot)@>[1]>>
$$
$$
\hbox{}\hskip -.4in\downarrow\hskip .7in\downarrow\hskip .6in\downarrow
$$
$$
\rightarrow R\Gamma_{{}_{\Theta^{-1}(D)}}(\Fdot)\rightarrow\Fdot\rightarrow R\Gamma_{{}_{\Theta^{-1}(\Bbb
R^2-D)}}(\Fdot)@>[1]>>
$$
proving $(\dagger)$ is equivalent to proving the isomorphism
$$H^*(R\Gamma_{{}_{\Theta^{-1}(\Bbb R^2-C)}}(\Fdot))_{\bold p}@>\cong>>H^*(R\Gamma_{{}_{\Theta^{-1}(\Bbb
R^2-D)}}(\Fdot))_{\bold p}.$$
Therefore, it suffices to produce a fundamental system $\{N_i\}$ of open neighborhoods of $\bold p$ such that we have
induced isomorphisms
$$\Bbb H^*(N_i-\Theta^{-1}(C);\ \Fdot)@>\cong>>\Bbb H^*(N_i-\Theta^{-1}(D);\ \Fdot). \tag{*}$$
Let $E:=\{u\geqslant 0, v\geqslant 0\}$. To show $(*)$ via moving the wall, we will produce $N_i$ such that $\Theta_{|_{N_i}}$
is a stratum-preserving, locally trivial fibration over $\Bbb R^2-C-E$; as $\Bbb R^2-C-E$ consists of two contractible
pieces, this will imply that the obvious homotopy from $\Bbb R^2-C$ to $\Bbb R^2-D$ lifts to give us $(*)$.

\vskip .2in

Take local embeddings of $X$ and $Y$ into affine spaces. Let $B_\epsilon$ denote a closed ball of radius $\epsilon$
centered at $\pi_1(\bold p)$ intersected with $X$, and let $B_\delta$ denote a closed ball of radius $\delta$
centered at $\pi_2(\bold p)$ intersected with $Y$. Let $\overset\circ\to B_\epsilon$ and $\overset\circ\to B_\delta$ denote
the intersections of the associated open balls with $X$ and $Y$, respectively. For positive $\eta$, let $T_\eta$ denote the
open square in
$\Bbb C$ given by $T_\eta =\{z\ |\ |\operatorname{Re}z|<\eta, \ |\operatorname{Im}z|<\eta\}$. We claim that
$$N_i:=
\big(\overset\circ\to B_{\epsilon_i}\times \overset\circ\to B_{\delta_i}\big)\cap (f\circ \pi_1)^{-1}(T_{\alpha_i})\cap
(g\circ \pi_2)^{-1}(T_{\beta_i}),$$
where $\alpha_i\ll\epsilon_i$ and $\beta_i\ll \delta_i$ is a fundamental system of open neighborhoods for which $(*)$ holds.

To see this, endow $B_{\epsilon_i}\times B_{\delta_i}$ with the obvious Whitney stratification, and consider the map $\Omega_i:
B_{\epsilon_i}\times B_{\delta_i}\rightarrow \Bbb C^2$ given by $\Omega_i:= (f\circ\pi_1, g\circ\pi_2)$. Let $a$ and $b$ denote
the coordinates in $\Bbb C^2$. The stratified critical points of $\Omega_i$ occur at points $(\bold x, \bold y)$ where either
$\bold x$ is in the stratified critical locus of
$f_{|_{B_{\epsilon_i}}}$ or $\bold y$ is in the stratified critical locus of
$g_{|_{B_{\delta_i}}}$; the standard Milnor fibration argument guarantees that, near $\bold p$, in a small neighborhood of
$\bold 0\in\Bbb C^2$, the stratified critical values occur only along $V(ab)$. Therefore, for $\alpha_i$ and $\beta_i$
sufficiently small, $\Omega_i$ is a proper stratified submersion over $T_{\alpha_i}\times T_{\beta_i}-V(ab)$. As
$\Theta_{|_{N_i}}=\operatorname{Re}(\Omega_i)$, it follows immediately that $\Theta_{|_{N_i}}$ is a stratum-preserving,
locally trivial fibration over $\Bbb R^2-C-E$.\qed

\vskip .3in

In the next proposition, we refer to $P$ and $Q$ and, of course, later we will apply this proposition to the $P$ and
$Q$ given in Lemma 1.2. However, in Proposition 1.3, $P$ and $Q$ are completely general.

\vskip .5in

\noindent{\bf Proposition 1.3}. {\it Let $r_1$ and $r_2$ denote the projections of $P\times Q$ onto $P$ and $Q$,
respectively. Let $l:P\times Q\rightarrow X\times Y$,
$i:P\rightarrow X$, and $j:Q\rightarrow Y$ be such that $i\circ r_1 = \pi_1\circ l$ and $i\circ r_2= \pi_2\circ l$, i.e., we
have a commutative diagram
$$ P@<r_1<< P\times Q@>r_2>> Q
$$
$$
\hbox{}\ \downarrow i\hskip .3in\downarrow l\hskip .4in\downarrow j
$$
$$ X@<\pi_1<< X\times Y@>\pi_2>> Y.
$$ 

Then, there is a natural isomorphism
$$ i^!\Adot\lboxtimes j^!\Bdot\ \cong\ l^!(\Adot\lboxtimes\Bdot).
$$

}

\vskip .2in

\noindent{\it Proof}. We use 3.4.4 and 3.1.13 of [{\bf K-S}].

$$ i^!\Adot\lboxtimes j^!\Bdot\ \cong\ \Cal Di^*\Cal D\Adot\lboxtimes j^!\Bdot\ \cong\ \rhom(r_1^*i^*\Cal D\Adot,
r_2^!j^!\Bdot)\ \cong
$$
$$
\rhom(l^*\pi_1^*\Cal D\Adot, l^!\pi_2^!\Bdot)\ \cong\ l^!\rhom(\pi_1^*\Cal D\Adot, \pi_2^!\Bdot)\ \cong\ l^!(\Cal D\Cal
D\Adot\lboxtimes \Bdot)\ \cong\ l^!(\Adot\lboxtimes\Bdot).\qed
$$

\vskip .3in

In the next proposition, we use repeatedly that if $t: T\hookrightarrow W$ is the inclusion of a closed subset and
$\Fdot$ is a complex on $W$, then
$R\Gamma_{{}_T}(\Fdot)$ is naturally isomorphic to $t_!t^!\Fdot$; this follows from 3.1.12 of [{\bf K-S}].

\vskip .5in

\noindent{\bf Proposition 1.4}. {\it We continue with the notation from the previous proposition, but we now assume that $l$,
$i$, and $j$ are inclusions of closed subsets. Then, there is a natural isomorphism
$$ R\Gamma_{{}_P}(\Adot)\lboxtimes R\Gamma_{{}_Q}(\Bdot)\ \cong\ R\Gamma_{{}_{P\times Q}}(\Adot\lboxtimes\Bdot).
$$

}

\vskip .2in

\noindent{\it Proof}. Let $\check\imath:P\times Y\hookrightarrow X\times Y$ denote the inclusion and let
$\check\pi_1:P\times Y\rightarrow P$ denote the projection. Analogously, let $\check\jmath:X\times Q\hookrightarrow X\times
Y$ denote the inclusion and let
$\check\pi_2:X\times Q\rightarrow Q$ denote the projection.

We have a natural isomorphism
$$ R\Gamma_{{}_P}(\Adot)\lboxtimes R\Gamma_{{}_Q}(\Bdot)\ \cong\ \pi_1^*i_!i^!\Adot\lotimes \pi_2^*j_!j^!\Bdot.
$$ Using the dual of 3.1.9 of [{\bf K-S}], $\pi_1^*i_!\cong \check\imath_!\check\pi_1^*$ and $\pi_2^*j_!\cong
\check\jmath_!\check\pi_2^*$. Thus, we have
$$ R\Gamma_{{}_P}(\Adot)\lboxtimes R\Gamma_{{}_Q}(\Bdot)\ \cong\ \check\imath_!\check\pi_1^*i^!\Adot\lotimes
\check\jmath_!\check\pi_2^*j^!\Bdot,
$$ and by the K\"unneth formula (Exercise II.18.i of [{\bf K-S}]), this last expression is naturally isomorphic to
$l_!(i^!\Adot\lboxtimes j^!\Bdot)$. Apply Proposition 1.3.\qed

\vskip .5in

\noindent{\bf 1.5 (Proof of the Sebastiani-Thom Isomorphism)}. We will use all of our previous notation and results. Let
$s_1:V(f)\hookrightarrow X$ and $s_2:V(g)\hookrightarrow Y$  denote the inclusions.

Then,
$$
\phi_{{}_f}\Adot\lboxtimes \phi_{{}_g}\Bdot \ = \ p_1^*\phi_{{}_f}\Adot\lotimes p_2^*\phi_{{}_g}\Bdot\ \cong\
p_1^*s_1^*R\Gamma_{{}_P}(\Adot)\lotimes p_2^*s_2^*R\Gamma_{{}_Q}(\Bdot)\ \cong
$$
$$ m^*\pi_1^*R\Gamma_{{}_P}(\Adot)\lotimes m^*\pi_2^*R\Gamma_{{}_Q}(\Bdot)\ \cong\ m^*(R\Gamma_{{}_P}(\Adot)\lboxtimes
R\Gamma_{{}_Q}(\Bdot))\ \cong
$$
$$ m^*R\Gamma_{{}_{P\times Q}}(\Adot\lboxtimes\Bdot) \ \cong\ k^*\phi_{{}_\map}\big(\Adot\lboxtimes\Bdot\big).
$$

\vskip .1in

The remaining statements of the theorem -- other than the monodromy statement follow from Lemma 1.1.

 The monodromy statement
follows at once from the proof  of the Sebastiani-Thom Isomorphism,  for the monodromy of $f\circ\pi_1+g\circ\pi_2$ results
from how the set $\{(\bold x, \bold y)\in X\times Y\ |\
 \operatorname{Re}(e^{-i\theta}(f(\bold x)+g(\bold y)))\leqslant 0\}$ moves as $\theta$ goes from $0$ to $2\pi$. The
isomorphism in Lemma 1.2 identifies this with the movement of $\{e^{-i\theta}\operatorname{Re}f(\bold x)\leqslant
0\}\times\{e^{-i\theta}\operatorname{Re}g(\bold y)\leqslant 0\}$, which describes the product of the two monodromies of
$f$ and $g$.\qed

\vskip .2in

\noindent \S2. {\bf Consequences}

\vskip .2in

Certainly there is some satisfaction in knowing that the Sebastiani-Thom Isomorphism holds for general spaces, even with
constant $\Bbb Z$-coefficients and only on the level of cohomology; in this case, the isomorphism reduces to
$$
\widetilde H^{i-1}(F_{{}_{\map, \bold p}})\cong 
$$
$$\Big[\ \bigoplus_{a+b=i}\Big(\widetilde H^{a-1}(F_{{}_{f, \pi_1(\bold
p)}})\otimes\widetilde H^{b-1}(F_{{}_{g, \pi_2(\bold p)}})\Big)\Big]\ \oplus\ \Big[\
\bigoplus_{c+d=i+1}\operatorname{Tor}\Big(\widetilde H^{c-1}(F_{{}_{f, \pi_1(\bold p)}}),\ \widetilde H^{d-1}(F_{{}_{g,
\pi_2(\bold p)}})\Big)\Big],
$$
where, as in the introduction, $\widetilde H$ denotes reduced, integral cohomology and
$F_{{}_{h, \bold x}}$ denotes the Milnor fibre of $h$ at $\bold x$.

\vskip .2in

Suppose, however, that we restrict ourselves to using coefficients in a field, say $\Bbb Q$ or $\Bbb C$. Then, when the
external tensor product $\lboxtimes$ is applied to two perverse sheaves, it returns a perverse sheaf. In addition, the
vanishing cycle functor $\phi_h$ takes perverse sheaves to perverse sheaves. Therefore, the Sebastiani-Thom Isomorphism yields
an isomorphism in the Abelian category of perverse sheaves on $X\times Y$ and, consequently the isomorphism preserves much
more subtle data than that provided by the stalk cohomology.

T. Braden works in this context in [{\bf B}], and it was his Lemma 3.16 which motivated the writing of this paper.   In [{\bf
B}], the base ring $R$ is the field $\Bbb C$,
$Y:=\Bbb C^m$,
$g:\Bbb C^m\rightarrow\Bbb C$ is the ordinary quadratic function $g:= y_1^2+\dots +y_m^2$, and $\Bdot:=\Bbb C^\bullet_Y$, the
constant sheaf on $Y$.  Hence, $\phi_g\Bdot$ is the constant sheaf on the origin, shifted by $-m$, and extended by zero onto
$V(g)$, i.e., if $\alpha$ denotes the inclusion of $\bold 0$ into $\Bbb C^m$, $\phi_g\Bdot\cong \alpha_!\Bbb
C^\bullet_{{}_{\bold 0}}[-m]$. Let $\tau:V(f)\hookrightarrow V(\map)$ be the inclusion given by $\tau(\bold x):=(\bold x,
\bold 0)$.

\vskip .2in

If $X$, $f$, and $\Adot$ are still arbitrary, we wish to show

\vskip .3in

\noindent{\bf Corollary}. {\it There is a natural isomorphism
$$
\phi_{{}_\map}(\pi_1^*\Adot)\cong \tau_*(\phi_{{}_f}\Adot)[-m].
$$

}

\vskip .2in

\noindent{\it Proof}. As the critical locus of $g$ is simply the origin, it follows from Lemma 1.1 and the Sebastiani-Thom
Isomorphism that

$$
\phi_{{}_\map}(\pi_1^*\Adot)\ \cong\ k_!(p_1^*\phi_{{}_f}\Adot\lotimes p_2^*\alpha_!\Bbb
C^\bullet_{{}_{\bold 0}}[-m]).
$$

Consider the pull-back diagram
$$
V(f)\times \bold 0@>\ \hat\alpha\ \ \ >>V(f)\times V(g)
$$
$$
\hat p_2\downarrow\hskip .9in\downarrow p_2
$$
$$
\hbox{}\hskip .2in\bold 0@>\hskip .3in\alpha\hskip .3in>>V(g).
$$
Then, we have $p_2^*\alpha_!\cong \hat\alpha_!\hat p_2^*$, and so 
$$
\phi_{{}_\map}(\pi_1^*\Adot)\ \cong\ k_!(p_1^*\phi_{{}_f}\Adot\lotimes \hat\alpha_!\hat p_2^*\Bbb
C^\bullet_{{}_{\bold 0}}[-m])\ \cong\ k_!(p_1^*\phi_{{}_f}\Adot\lotimes \hat\alpha_!\Bbb
C^\bullet_{{}_{V(f)\times\bold 0}}[-m]).
$$
Applying Proposition 2.6.6 of [{\bf K-S}], we obtain that
$$
\phi_{{}_\map}(\pi_1^*\Adot)\ \cong\ k_!(p_1^*\phi_{{}_f}\Adot\lotimes \hat\alpha_!\Bbb
C^\bullet_{{}_{V(f)\times\bold 0}}[-m])\ \cong
$$
$$k_!\hat\alpha_!(\hat\alpha^* p_1^*\phi_{{}_f}\Adot\lotimes\Bbb
C^\bullet_{{}_{V(f)\times\bold 0}}[-m])\ \cong\ k_!\hat\alpha_!(\hat\alpha^* p_1^*\phi_{{}_f}\Adot[-m])\ \cong\
(k\circ\hat\alpha)_!(p_1\circ\hat\alpha)^*\phi_{{}_f}\Adot[-m].
$$
As $p_1\circ\hat\alpha$ is the isomorphism which identifies $V(f)\times\bold 0$ and $V(f)$, and as $k\circ\hat\alpha$ is the
closed inclusion of $V(f)\times\bold 0$ into $V(\map)$, we are finished.\qed

\vskip .4in

We thank Tom Braden for motivating this work and for many useful comments during the preparation of this paper.

\enddocument